\documentclass[12pt,a4paper]{article}
\usepackage{amsmath,amssymb,amsfonts}
\def\Bbb{\mathbb}

\title{\bf Equality of Dedekind sums modulo $24\mathbb Z$}

\author{Kurt Girstmair}
\date{}

\makeatletter
\let\@@maketitle=\maketitle
\def\maketitle{\def\thispagestyle##1{\relax}\@@maketitle}
\makeatother
\newtheorem{theorem}{Theorem}
\newtheorem{prop}{Proposition}

\def\BE{\begin{equation}}
\def\EE{\end{equation}}
\def\BD{\begin{displaymath}}
\def\ED{\end{displaymath}}
\def\BA{\begin{array}}
\def\EA{\end{array}}
\def\BEA{\begin{eqnarray*}}
\def\EEA{\end{eqnarray*}}
\def\BI{\bibitem}

\def\N{\Bbb N}
\def\Z{\Bbb Z}

\def\R{\Bbb R}

\def\phi{\varphi}
\def\EPS{\varepsilon}

\def\MOD{\: \mbox{mod} \:}

\def\MB{\mbox}
\def\LD{\ldots}

\def\DIV{\,|\,}
\def\NDIV{\, \nmid \,}

\def\MN{\medskip\noindent}
\def\STOP{\hfill$\Box$}

\def\JS#1#2{ \left( \frac{#1}{#2} \right) }
\def\DED{Dedekind }

\begin{document}
\maketitle

\begin{abstract}

\noindent
Let $S(a,b)=12s(a,b)$, where $s(a,b)$ denotes the classical \DED sum.
In a recent note E. Tsukerman gave a necessary and sufficient condition for $S(a_1,b)-S(a_2,b)\in 8\Z$.
In the present paper we show that this condition is equivalent to $S(a_1,b)-S(a_2,b)\in 24\Z$, provided that $9\NDIV b$.
Tsukerman also obtained a congruence mod 8 for $bT(a,b)$, where $T(a,b)$ is the alternating sum of the
partial quotients of the continued fraction expansion of $a/b$. We show that the respective congruence holds mod
$24$ if $3\NDIV b$ and mod $72$ if $3\DIV b$.
\end{abstract}

\section*{1. Introduction and results}

Let $a$ be an integer, $b$ a natural number, and $(a,b)=1$. The classical \DED sum $s(a,b)$ is defined by
\BD
   s(a,b)=\sum_{k=1}^{b} ((k/b))((ak/b)).
\ED
Here
\BD
  ((x))=\begin{cases}
                 x-\lfloor x\rfloor-1/2 & \MB{ if } x\in\R\smallsetminus \Z; \\
                 0 & \MB{ if } x\in \Z
        \end{cases}
\ED
(see \cite[p. 1]{RaGr}).
It is often more
convenient to work with
\BD
 S(a,b)=12s(a,b)
\ED instead (see, for instance, formula (\ref{2.6}) below). Since $S(a+b,b)=S(a,b)$, we obtain all \DED sums if $a$ is restricted to the
range $0\le a\le b-1$, $(a,b)=1$.

In the recent note \cite{Ts}, E. Tsukerman gave a necessary and sufficient condition for the
equality of $S(a_1,b)$ and $S(a_2,b)$ modulo $8\Z$. This condition involves the function
$\mu$, which is defined, for $a,b$ as above, as follows
\BD
\label{1.2}
\mu(a,b)=\begin{cases}
           2-2\JS ab, & \MB{ if } b \MB{ is odd;}\\
           (a-1)(a+b-1), & \MB{ if } b \MB{ is even.}
         \end{cases}
\ED
Here $\JS ab$ is the Jacobi symbol.
Tsukerman's condition is phrased by means of the residue class
\BD
\label{1.3}
b(a_2\mu(b, a_1)-a_1\mu(b,a_2)) \enspace \MB{ mod } 8b.
\ED
We observe, however, that this residue class depends only on the residue classes of $\mu(b,a_1)$ and $\mu(b,a_2)$ modulo 8,
not of the values of $\mu(b,a_1)$ and $\mu(b,a_2)$ themselves.
Therefore, we may replace the function $\mu$ by the following simpler function, which we henceforth also call $\mu$.
\BD
\label{1.6}
\mu(a,b)=\begin{cases}
           2-2\JS ab, & \MB{ if } b \MB{ is odd;}\\
           4, & \MB{ if } b \equiv 0\mod 4 \MB{ and } a\equiv 3\mod 4;\\
           0, & \MB{ otherwise.}
         \end{cases}
\ED
In this paper we show

\begin{theorem} 
\label{t1}
Let $a_1, a_2\in\N$ be relatively prime to $b\in\N$. Suppose, further, that $9\NDIV b$.
Then
\BD
  S(a_1,b)-S(a_2,b)\in 24\Z
\ED
if, and only if,
\BE
\label{1.4}
b(a_2\mu(b, a_1)-a_1\mu(b,a_2))\equiv (a_1-a_2)(b-1)(a_1a_2+b-1)\mod 8b.
\EE
\end{theorem} 

This equivalence cannot be extended to the case $9\DIV b$ in an obvious way, as we show in Section 3.
Tsukerman showed that (\ref{1.4}) is equivalent to $S(a_1,b)-S(a_2,b)\in 8\Z$ for arbitrary natural numbers $b$,
i.e., he needed not assume $9\NDIV b$.

For $a\in \Z$ and $b\in \N$, let
\BD
  \frac ab=[a_0,a_1,\LD,a_n]
\ED
denote the regular continued fraction expansion of $a/b$. The partial quotients $a_1,\LD,a_n$ are natural numbers. We do not assume $a_n\ge 2$,
but require $n$ to be odd, instead. Define
\BE
\label{1.8}
   T(a,b)=\sum_{k=0}^n(-1)^{k-1}a_k.
\EE
In the said paper, Tsukerman showed, for $a,b\in\N$, $(a,b)=1$,
\BE
\label{1.10}
 bT(a,b)\equiv -\mu(a,b)+b^2+2-a-a^*\mod 8,
\EE
with $a^*\in\{1,\LD,b-1\}$, $aa^*\equiv 1$ mod $b$. Our new definition of $\mu$ suggests a more explicit form of (\ref{1.10}),
which we use in the following Theorem. To this end we define $\EPS\in\{\pm 1\}$ by the congruence
\BE
\label{1.12}
   a\equiv\EPS \mod 3
\EE
for each $a\in \Z$, $3\NDIV a$.

\begin{theorem} 
\label{t2}

Let $a\in\Z$ be relatively prime to $b\in\N$.\\
{\rm (a)} Let $b$ be odd. If $3\NDIV b$, then
\BD
  bT(a,b)\equiv 9+18 \JS ab -a -a^* \mod 24.
\ED
If $3\DIV b$, then
\BD
  bT(a,b)\equiv 9+18 \JS ab-16\EPS -a -a^* \mod 72.
\ED
{\rm (b)} Let $b\equiv 2$ mod $4$ or let both $b\equiv 0$ mod $4$ and $a\equiv 3$ mod $4$ hold.
If $3\NDIV b$, then
\BD
  bT(a,b)\equiv 6 -a -a^* \mod 24.
\ED
If $3\DIV b$, then
\BD
  bT(a,b)\equiv 54-16\EPS -a -a^* \mod 72.
\ED
{\rm (c)} Let $b\equiv 0$ mod $4$ and $a\equiv 1$ mod $4$. If $3\NDIV b$, then
\BD
  bT(a,b)\equiv 18 -a -a^* \mod 24.
\ED
If $3\DIV b$, then
\BD
  bT(a,b)\equiv 18-16\EPS -a -a^* \mod 72.
\ED

\end{theorem} 

In Section 3 we exhibit many examples that illustrate both Theorem \ref{t1} and the fact that this theorem does
not hold if $9\DIV b$.

\section*{2. Proofs}

Our main tools are two congruences modulo 3 for \DED sums. First we observe that $bS(a,b)$ is an integer;
moreover, if $3$ does not divide $b$, then
\BE
\label{2.2}
  bS(a,b)\equiv 0\mod 3.
\EE
These  assertions follow from \cite[p. 27, Th. 2]{RaGr}). On the other hand, if $3\DIV b$,
\BE
\label{2.4}
 bS(a,b)\equiv 2\EPS\mod 9,
\EE
where $\EPS$ is defined as in (\ref{1.12}) (see \cite[formula (70)]{Sa}).

\MN
{\em Proof of Theorem \ref{t1}.}
Suppose, first, that $3\NDIV b$.
Because of (\ref{2.2}), we may write
\BD
  S(a_1,b)=\frac{3k_1}b,\enspace S(a_2,b)=\frac{3k_2}b
\ED
with integers $k_1$, $k_2$.  By \cite[Th. 3.1]{Ts}, the congruence (\ref{1.4}) is equivalent to $S(a_1,b)-S(a_2,b)\in 8\Z$.
Accordingly, (\ref{1.4}) is also equivalent to
\BD
  \frac{3(k_1-k_2)}b=8r, \enspace r\in \Z.
\ED
However, $3\NDIV b$, and so this means $3\DIV r$. This proves Theorem \ref{t1} in the case $3\NDIV b$.
Suppose now that $3\DIV b$. Then the congruence (\ref{1.4}) implies
$(a_1-a_2)(a_1a_2-1)\equiv 0\MOD 3$. Hence we obtain, from (\ref{2.4})
\BD
  S(a_1,b)=\frac{2\EPS+9k_1}b,\enspace S(a_2,b)=\frac{2\EPS+9k_2}b
\ED
with a common value $\EPS\equiv a_1\equiv a_2$ mod 3 and $k_1,k_2\in\Z$.
Accordingly, (\ref{1.4}) is equivalent to
\BD
  \frac{9(k_1-k_2)}b=8r, \enspace r\in\Z.
\ED
If $9\DIV b$, this simply means $S(a_1,b)-S(a_2,b)\in 8\Z$, so this is just Tsukerman's result.
However, if $9\NDIV b$, we obtain $3\DIV r$, which yields the theorem in the case $3\DIV b, 9\NDIV b$.
\STOP

\MN
{\em Proof of Theorem \ref{t2}.}
The Barkan-Hickerson-Knuth formula says
\BE
\label{2.6}
  S(a,b)=T(a,b)+\frac{a+a^*}b-3
\EE
(see, for instance, \cite{Hi}).
Note that this formula is often enunciated only for the case $0\le a<b$, but it is, in fact, valid
for arbitrary integers $a$ relatively prime to $b$, provided that $T(a,b)$ is defined as in (\ref{1.8}).
Hence we obtain, by (\ref{2.2}) and (\ref{2.6}),
\BE
\label{2.8}
  bT(a,b)\equiv -a-a^*\mod 3
\EE
if $3\NDIV b$. In the case $3\DIV b$, (\ref{2.4}) and (\ref{2.6}) give
\BE
\label{2.10}
  bT(a,b)\equiv 2\EPS-a-a^*\mod 9
\EE
instead. Further, Tsukerman's congruence (\ref{1.10})
is also valid for arbitrary integers $a$ relatively prime to $b$, as we easily check.
We combine (\ref{1.10}) with the congruences (\ref{2.8}), (\ref{2.10}) by means of the Chinese remainder theorem.
This readily gives Theorem \ref{t2}. \STOP

\section*{3. A proposition yielding examples}

Our examples arise from the following proposition.

\begin{prop} 
\label{p1}
Let $c,d$ be odd natural numbers, $d\ge 3$. Put $b=cd^2$ and $a=cd+1$. Then
\BD
  S(1,b)-S(a,b)= c(d^2-1).
\ED
\end{prop} 

\MN
{\em Proof.} We apply the reciprocity law for \DED sums (see \cite[p. 5]{RaGr}), which gives
\BD
 S(a,b)=-S(b,a)+ \frac{a}b+\frac b{a}+\frac 1{ab}-3.
\ED
Now $b\equiv-d\mod a$, hence the reciprocity law says
\BD
S(b,a)=S(-d,a)=S(a,d)-\frac{a}d-\frac d{a}-\frac 1{ad}+3.
\ED
However, $a\equiv 1\mod d$, and so $S(a,d)=S(1,d)=d-3+2/d$. Inserting the values $b=cd^2$ and $a=cd+1$ gives
\BD
  S(a,b)=c-3+\frac2{b}.
\ED
Since $S(1,b)=b-3 +2/b$, we obtain the desired result.
\STOP

\medskip
In the setting of the proposition, let $3\NDIV d$. Then $d^2-1\equiv 0\mod 24$, so the proposition yields many examples
with $S(1,b)-S(a,b)>0$ and $S(1,b)-S(a,b)\equiv 0$ mod 24. On the other hand, if $3\DIV d$, then $d^2-1\equiv 0$ mod 8,
but $d^2-1\not\equiv 0$ mod 24. If, therefore, $3\NDIV c$, we obtain many examples with $S(1,b)-S(a,b)\equiv 0$ mod 8,
but $S(1,b)-S(a,b)\not\equiv 0$ mod 24.


\vspace{0.5cm}
\noindent
Kurt Girstmair            \\
Institut f\"ur Mathematik \\
Universit\"at Innsbruck   \\
Technikerstr. 13/7        \\
A-6020 Innsbruck, Austria \\
Kurt.Girstmair@uibk.ac.at

\end{document}